%% file: agt-2-19.tex
\documentclass{gtart}

\input agtout

\lognumber{19}
\volumenumber{2}
\volumeyear{2002}
\papernumber{19}
\published{22 May 2002}
\pagenumbers{391}{402}
\received{19 February 2002}
\accepted{12 May 2002}

\usepackage{amsmath, amssymb}

\let\Bbb\mathbb
\theoremstyle{plain}
\newtheorem{thm}{Theorem}[section]
\newtheorem{lem}[thm]{Lemma}

\newtheorem{prop}[thm]{Proposition}

\theoremstyle{definition}

\def \CPb {\overline{\mathbf{CP}}^{\,2}}
\def \CP {{\mathbf{CP}}^{2}}

\def \Z {\mathbb{Z}}
\def \Sig{\Sigma}

\def \SS {\mathbf{S}^2\times \mathbf{S}^2}

\def \a {\alpha}
\def \b {\beta}
\def \g {\gamma}

\def \L {\Lambda}

\def \s {\sigma}

\def \z {\zeta}

\def \bd {\partial}
\def \x {\times}
\def \- {\setminus}
\def \C {\subset}
\def \ve {\varepsilon}

\def \ch {\chi_{{}_h}}
\def \lr {\leftrightarrow}

\def \sign {\text{sign}}

\begin{document}

\title{Rational surfaces and symplectic $4$-manifolds\\with one basic class}

\authors{Ronald Fintushel\\Jongil Park\\Ronald J. Stern}
\shortauthors{Fintushel, Park and Stern}

\address{Department of Math, Michigan State University,
East Lansing, MI 48824, USA\\\medskip\\Department of Math,  
Konkuk University,
1 Hwayang-dong\\Kwangjin-gu, Seoul 143-701, Korea\\\medskip\\Department 
of Math, University of California,
Irvine, CA 92697, USA}

\asciiaddress{Department of Mathematics, Michigan State University\\
East Lansing, MI 48824, USA\\Department of Mathematics,  
Konkuk University \\
1 Hwayang-dong, Kwangjin-gu, Seoul 143-701, Korea\\Department 
of Mathematics, University of California \\Irvine,
CA 92697, USA}

\email{ronfint@math.msu.edu, jipark@konkuk.ac.kr, rstern@math.uci.edu}

\begin{abstract}  We present constructions of simply connected symplectic
4-manifolds which have (up to sign) one basic class and which fill up the
geographical region between the half-Noether and Noether lines.
\end{abstract}

\primaryclass{57R57}
\secondaryclass{57R17}

\keywords{$4$-manifold, Seiberg-Witten invariant}

\maketitlepage

\section{Introduction\label{Intro}}

For minimal complex surfaces $S$ of general type, the Noether
inequality states that
$c_1^2(S)\ge 2\ch(S)-6$, where $\ch(X)$ denotes the holomorphic Euler
number of $X$.
($\ch(X)=\frac{1}{4}(e(X)+\sign(X))$ where $e$ is the Euler characteristic and
$\sign$ is the signature of the intersection form.) The line
$c_1^2=2\ch-6$ in the
$(\ch, c_1^2)$-plane is often called the Noether line. In terms of
gauge theory, one
of most notable features of a minimal surface of general type is
that, up to sign, it
has exactly one (Seiberg-Witten) basic class \cite{W}. In \cite{rat} the first
and third authors produced examples of symplectic (see \cite{Sym})
4-manifolds with
one basic class which lie on the `half-Noether' line $c_1^2=\ch-3$. The
inability to construct examples (even smoothly) of 4-manifolds with
one basic class
and
$c_1^2<\ch-3$ led them to conjecture that such manifolds fail to
exist. Interest in
this problem was reignited recently by a paper of Mari\~{n}o, Moore,
and Peradze
\cite{MMP} which gave a plausibility argument via physics.

In the current article, we show the existence of symplectic manifolds
with one basic
class which fill the region in the $(\ch, c_1^2)$-plane between the
half-Noether
and Noether lines. Specifically we prove:

\begin{thm}\label{T} For every pair of positive integers $(x,c)$ with
$0< x-3 \le c
\le 2x-6$ there is a simply connected symplectic 4-manifold $X$ with
$c_1^2(X)=c$,
$\chi(X)=x$ and (up to sign) one basic class.
\end{thm}

We were drawn to this problem by a question of Paul Feehan. The
manifolds produced
in Theorem~\ref{T} serve to simplify some of the calculations necessary in the
Feehan-Leness program to prove the equivalence of the Seiberg-Witten
and Donaldson
invariants. Another goal of this paper is to exhibit further techniques for
constructing
4-manifolds. In light of this we will present two different 
constructions for the manifold in question. Each construction brings 
to light interesting properties of rational and elliptic surfaces.

The key to our constructions is to understand configurations of
embedded surfaces
in rational surfaces which can be rationally blown down. We close
this introduction
by reminding the reader of the notion of rational blowdown. (See \cite{rat} and
\cite{P}.) Let
$C_n$ denote the simply connected smooth
$4$-manifold with boundary obtained by plumbing $n-1$ disk bundles over the
$2$-sphere according to the linear diagram:

\vspace{-4mm}
\centerline{\unitlength .9cm
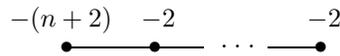
\begin{picture}(5,2)\small
\put(.9,.7){$\bullet$}
\put(1,.8){\line(1,0){1.3}}
\put(2.2,.7){$\bullet$}
\put(2.3,.8){\line(1,0){.75}}
\put(3.3,.8){.}
\put(3.5,.8){.}
\put(3.7,.8){.}
\put(4,.8){\line(1,0){.75}}
\put(4.65,.7){$\bullet$}
\put(.15,1.1){$-(n+2)$}
\put(2.1,1.1){$-2$}
\put(4.55,1.1){$-2$}
\end{picture}}
\vspace{-4mm}

Here, each node denotes a disk bundle over $S^2$  with Euler
class indicated by the label; an interval indicates that the endpoint
disk bundles
are plumbed,  i.e. identified fiber to base.
Label the homology classes represented by the spheres in $C_n$ by
$S_0,S_1,\dots,S_{n-2}$ so that the self-intersections are
$S_0^2=-(n+2)$ and, for $j=1,\dots,n-2$, $S_j^2=-2$. Further, orient
the spheres so that $S_j\cdot S_{j+1}=+1$.  Then $C_n$ is a  $4$-manifold with
negative definite intersection form and with boundary the lens space
$L(n^2,1-n)$.
The lens space $L(n^2,1-n)=\bd C_n$ bounds a rational ball $B_n$ with \
$\pi_1(B_n)={\Z}_n$ and a surjective inclusion-induced homomorphism \
$\pi_1(L(n^2,1-n)={{\Z}}_{n^2}\to \pi_1(B_n)$. If $X$ is a smooth $4$-manifold
containing an embedded copy of $C_n$, its `rational blowdown' is the result of
replacing $C_n$ by the rational ball $B_n$. Rationally blowing down
$C_n$ increases
$c_1^2$ by $n-1$ but does not change $\ch$. It is a theorem of
Symington \cite{Sym}
that if the ambient manifold is symplectic, and each sphere $S_i$ of
a configuration
$C_n$ is a symplectic submanifold, the resultant manifold of the rational
blowdown is also symplectic.

A key result is:

\begin{thm}{\rm\cite{rat}}\label{taut}\qua Let $X$ be a simply connected
smooth 4-manifold with $b^+>1$ and containing the configuration
$C_n$. Suppose that all the (Seiberg-Witten) basic classes $k$ of $X$
satisfy
\[ k\cdot S_i=0, \ i=1,\dots, n-2, \ \ {\text{and}} \ \ |k\cdot S_0|\le n.\]
Then the result of rationally blowing down $C_n$ is a smooth
4-manifold whose basic
classes are in one-to-one correspondence with the basic classes $k$ of $X$
satisfying $|k\cdot S_0| = n$.
\end{thm}

The authors gratefully acknowledge the following external support. The first
author was partially supported NSF Grant DMS0072212, the second author by KOSEF
981-0103-015-2, and the third author by NSF Grant DMS9971667.

\section {Line arrangements and 4-manifolds}

In this section we shall construct rational surfaces
which contain Riemann surfaces of self-intersection $0$, along which
one is able
to form fiber sums. The result of these fiber sums will be elliptic
surfaces, Horikawa surfaces, and symplectic manifolds which sit on
the ``half-Noether
line'' $c_1^2=\ch - 3$. There are certainly other constructions of
these manifolds
(cf.\ \cite{rat}) and we shall describe one such in the next section,
however the
description below is the most convenient for our purposes.

Let $q$ be an integer $\ge 4$. To construct the first rational
surface, consider the
arrangement of $q$ lines in $\CP$ formed by taking $q-2$ lines through a common
point and two more lines in general position. Blow up the multiple point $x_0$
to obtain a configuration of rational curves in $\CP\#\CPb$ representing
$qH-(q-2)E$, where $H$ denotes the class of a line and $E$ the exceptional
curve. Smooth double points to obtain a smooth embedded holomorphic curve of
self-intersection $4q-4$ and genus  $q-2$ (as seen via the adjunction formula).
Now blow up $4q-4$ more points along the embedded surface to obtain
the rational
surface $R(q)$ with $c_1^2(R(q))=12-4q$ and with a surface
$\Sig_{R(q)}$ of genus
$q-2$ with trivial normal bundle. Furthermore, since an exceptional
curve $E_i$,
($i=1,\dots, 4q-4$) is a 2-sphere that intersects $\Sig_{R(q)}$ in
one point, the
complement, $R(q)\- \Sig_{R(q)}$ is simply connected.

To construct the second rational surface, start with the arrangement
of lines in
$\CP$ obtained by taking $p-3$ lines in $\CP$ meeting in one point
and then adding
three more lines in general position. Blow up the point of multiplicity $p-3$
to obtain a configuration of rational curves in $\CP\#\CPb$ representing
$pH+(p-3)E$. After smoothing the double points, we obtain a smooth embedded
holomorphic curve of self-intersection $6p-9$ and genus $2p-5$.
Finally, blow up
$6p-9$ points along the embedded surface to obtain the rational
surface $S(p)$ with
$c_1^2(S(p))=17-6p$ and with a surface $\Sig_{S(p)}$ of genus $2p-5$
with trivial
normal bundle. Note also that $\Sig_{S(p)}$ intersects all the
exceptional classes,
and $S(p)\- \Sig_{S(p)}$ is simply connected.

Define $X_p$ to be the symplectic 4-manifold obtained by taking the
fiber sum of
$R(2p-3)$ and $S(p)$ along $\Sig_{R(2p-3)}$ and $\Sig_{S(p)}$. (Note that both
these surfaces have genus $2p-5$.) For fiber sums along surfaces of
genus $g$, one
has the general formulas

\cl{$c_1^2(A\#_{\Sig}B)= c_1^2(A) + c_1^2(B) + (8g-8)$,}
\vspace{1mm}
\cl{$\ch(A\#_{\Sig}B)= \ch(A) + \ch(B) + (g-1)$.}
\vspace{2mm}

It follows that $c_1^2(X_p)= 2p-7$ and $\ch(X_p)=2p-4$;
so $c_1^2(X_p)=\ch(X_p)-3$. Since the complements of $\Sig_{R(2p-3)}$ and
$\Sig_{S(p)}$ in $R(2p-3)$ and $S(p)$ are simply connected, so is $X_p$.

These manifolds, $X_p$ all have holomorphic Euler number $\ch(X_p)$ even.  To
obtain examples with odd $\ch$, modify the above construction as follows:
Start once more with the arrangement consisting of $p-3$ lines through a single
point and $3$ further lines in general position. Blow up the multiple point of
multiplicity $p-3$ and also one of the double points to obtain a
configuration of
rational curves in $\CP\#2\,\CPb$ representing the homology class
$pH-(p-3)E-2E_1$.
After smoothing the double points of the configuration one obtains a
smooth embedded
holomorphic curve of self-intersection $6p-13$ and genus $2p-6$. Blow
up $6p-13$
points along the embedded surface to obtain the rational surface $S'(p)$ with
$c_1^2(S'(p))=20-6p$ and with a surface $\Sig_{S'(p)}$ of genus
$2p-6$ which has a
trivial normal bundle.

Define $X'_p$ to be the symplectic 4-manifold obtained by taking the
fiber sum of
$R(2p-4)$ and $S'(p)$ along the genus $2p-6$ surfaces $\Sig_{R(2p-4)}$ and
$\Sig_{S'(p)}$. Then
$c_1^2(X'_p) = 2p-8$ and $\ch(X_p)=2p-5$;
so again, $c_1^2(X'_p)=\ch(X'_p)-3$, and as above, $X'_p$ is simply connected.

\section{Construction via rational blowdowns}

In order to compute the Seiberg-Witten invariants of the symplectic
4-manifolds $X_p$
and $X'_p$, it is useful to have an alternative construction. We
first concentrate
on $X_p$. Let $R=R(2p-3)$ and $\Sig_R=\Sig_{R(2p-3)}$, and let $S=S(p)$ and
$\Sig_S=\Sig_{S(p)}$. Then $\Sig_R$ represents the homology class
\[ (2p-3)H-(2p-5)E-\sum_{i=1}^{8p-16}E_i \ \in \ H_2(\CP\#\CPb\#(8p-16)\CPb).\]

The rational surface $R$ contains
the configuration $C=C_{2p-6}$ which is a linear plumbing of $2p-7$ holomorphic
spheres:

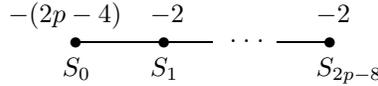
\centerline{\unitlength .9cm\small
\begin{picture}(5,2)
\put(.9,.7){$\bullet$}
\put(1,.8){\line(1,0){1.3}}
\put(2.2,.7){$\bullet$}
\put(2.3,.8){\line(1,0){.75}}
\put(3.3,.8){.}
\put(3.5,.8){.}
\put(3.7,.8){.}
\put(4,.8){\line(1,0){.75}}
\put(4.65,.7){$\bullet$}
\put(0,1.1){$-(2p-4)$}
\put(.8,.3){$S_0$}
\put(2.1,.3){$S_1$}
\put(4.55,.3){$S_{2p-8}$}
\put(2.1,1.1){$-2$}
\put(4.55,1.1){$-2$}
\end{picture}}
\vspace{-4mm}

where
\[ S_0=H-\sum_{i=1}^{2p-3}E_i,\ \ S_1=E_{2p-3}-E_{2p-2},\ \dots,\ \
S_{2p-8}=E_{4p-12}-E_{4p-11}. \]
Notice that $\Sig_R$ is disjoint
from the configuration $C$. This configuration can be
rationally blown down by replacing it with a rational ball $B_{2p-6}$ with
$\pi_1=\Z_{2p-6}$.  We claim that the rational surface
$S$ is the result of rationally blowing down $C$. Since $\Sig_R$ is
contained in the
complement of $C$, it gives rise to a surface in the new manifold.

\begin{prop}\label{P} Rational blowdown of the configuration $C$ in
$R$ yields $S$,
and the surface $\Sig_R$ becomes $\Sig_S \C S$.
\end{prop}

\begin{proof}
We shall prove this by rationally blowing down $C$ together with
$6p-9$ exceptional
curves in $R$. The result will be $\CP\#\CPb$, and $\Sig_R$ will get
blown down to
the class $p\, h - (p-3)\, e$, where $h$ and $e$ represent the
obvious classes in
$\CP\#\CPb$. (We shall use lower case notation in order not to
confuse these classes
with those used in the description of $\Sig_R$.)

The $6p-9$ exceptional curves in $R$ to be blown down are
$\{ E_{4p-10},\dots,E_{8p-16}\}$, and
$\{H-E-E_1,\dots,H-E-E_{2p-4}\}$. These curves
are all disjoint from $C$ (and from each other). Thus, if we choose,
we may first
blow down all the exceptional curves and then rationally blow down the
configuration $C$. Blowing down the $E_j$, $j=4p-10,\dots, 8p-16$ we obtain
$\CP\#\CPb\#(4p-9)\CPb$ containing the blown down surface $\Sig_R'$
which represents
the homology class
$(2p-3)H-(2p-5)E-\sum_{i=1}^{4p-11}E_i$.

Next blow down the exceptional curves $H-E-E_i$, $i=1,\dots, 2p-4$.
The result is a
rational surface $Q$ which has
\[ \{\a= H-E,\, \b=H-\sum_{i=1}^{2p-4}E_i,\, E_{2p-3},\dots, E_{4p-11} \} \]
as a basis for $H_2(Q)$. (Here we have compacted notation. If we
denote the blow
down map $R\to Q$ by $\pi$, then we should write $\a=\pi_*(H-E)$, etc. This
abbreviated notation should not cause any confusion, and we will
continue to use it
below.)
Note that both $\a$ and $\b$ are represented by
holomorphic spheres. Furthermore, $\a$ has self-intersection 0, and
it intersects
$\b$ once, hence for any nonnegative integer $k$, the class $k\a+\b$ is also
represented by an embedded 2-sphere.
This series of blowdowns takes $\Sig_R'$ to a surface
$\Sig_R''$ which represents the class
\[(4p-7)H-(4p-9)E-2\,\sum_{i=1}^{2p-4}E_i-\sum_{j=2p-3}^{4p-11}E_j \  = \
(4p-9)\,\a +2\,\b-\ve\] in $H_2(Q)$ (where $\ve = \sum_{j=2p-3}^{4p-11}E_j$) .

In $Q$, the
sphere $S_0$ of the configuration $C$ is given by $S_0=\b-E_{2p-3}$.  The
configuration defines a subspace of the second homology whose
orthogonal complement
$H_2(C)^{\perp}$ has basis
$\{\g_1,\g_2\}$ where
\[ \g_1= (2p-5)\a+\b\ \  \text{and}\ \  \g_2= \a-\ve \] with
intersection matrix:
\[ \begin{pmatrix} 2p-5 & 1 \\ 1 & -(2p-7)
\end{pmatrix}\] Both $\g_1$ and $\g_2$ are represented by embedded holomorphic
2-spheres in $Q\- C$. We have already seen this for $\g_1$, and it is clear for
$\g_2$. Because $H_2(C)$ is negative definite, it follows easily that
$H_2(C)^{\perp}= H_2(Q\- C)$, and in terms of our generators,
$\Sig_R''= 2\g_1+\g_2$.

Rationally blow down $C$, replacing it with the rational ball
$B_{2p-6}$. The result
is a symplectic (\cite{Sym}) 4-manifold $Y=(Q\- C)\cup B_{2p-6}$, and
the classes
$\g_1$ and $\g_2$ rationally generate $H_2(Y)$. Since $\g_1$ is
represented by a
symplectic 2-sphere of self-intersection $2p-5>0$, it follows from a theorem of
McDuff \cite{M} that $Y$ must be $\CP\#\CPb$.

If we view $Y$ as the ruled surface $\Bbb{F}_{2p-7}$ with fiber class $f$ and
positive and negative section classes $s_+$ and $s_-$, then $\g_1$
and $\g_2$ are
identified as $\g_1=s_++f$ and $\g_2=s_-$. Note that this agrees with the model
presented in \cite{rat} where it is shown that $\Bbb{F}_{2p-7}$ is the union of
$B_{2p-6}$ and a regular neighborhood of spheres representing $s_++f$
and $s_-$.
Since in $\Bbb{F}_{2p-7}$ we have
\[ s_++f = (p-2)\,h-(p-3)\,e \ \ \text{and} \ \ s_- = (4-p)\,h+(p-3)\,e. \]
It follows that
\[ h = \frac12\,(\g_1+\g_2) \ \ \text{and} \ \ e = \frac{1}{2p-6}
((p-4)\,\g_1 +
(p-2)\,\g_2).\]

In this process, the surface $\Sig_R$ has been blown down to a genus
$2p-5$ surface
representing $2\g_1+\g_2= p\,h+(p-3)\,e$; so when we blow up $6p-9$ 
times, we get
$\Sig_S$, and this proves the proposition.
\end{proof}

Similarly, let $R'=R(2p-4)$, $S'=S'(p)$, and $\Sig'_R = \Sig_{R(2p-4)}$,
$\Sig'_S=\Sig_{S'(p)}$. Then $\Sig'_R$ represents the homology class
\[ (2p-4)H-(2p-6)E-\sum_{i=1}^{8p-20}E_i \ \in \ H_2(\CP\#\CPb\#(8p-20)\CPb)\]
and $R'$ contains the configuration $C'=C_{2p-7}$ composed of
\[ S'_0=H-\sum_{i=1}^{2p-4}E_i,\ \ S'_1=E_{2p-4}-E_{2p-3},\ \dots,\ \
S'_{2p-9}=E_{4p-14}-E_{4p-13}. \]

\begin{prop}\label{P'} Rational blowdown of the configuration $C'$ in
$R'$ yields
$S'$ and the surface $\Sig'_R$ becomes $\Sig'_S \C S'$.
\end{prop}
\begin{proof} This can be proved in a fashion similar to the proposition above.
After blowing down $E_{4p-12},\dots, E_{8p-20}$ and $H-E-E_2,\dots,
H-E-E_{2p-5}$,
all of which are orthogonal to $C'$, we are left with
$U=\CP\#\CPb\#(2p-7)\CPb$. A
basis for $H_2(U)$ is given by
\[ \{\a=H-E,\ \b=H-\sum_{i=2}^{2p-5}E_i,\ E_1,\ E_{2p-4},\dots, E_{4p-13}\}, \]
and $H_2(C')^{\perp}$ is generated by $\{\z=\a-E_1,\ \a-\ve,\
(2p-7)\a+\b\}$, where
$\ve=\sum_{2p-4}^{4p-13}E_j$. The surface $\Sig'_R$ gets blown down
to a surface
$\Sig'_{R,U}$ which, in terms of this basis, represents the class
$\z+(\a-\ve)+2((2p-7)\a+\b)$.

Rationally blow down $C'$ to obtain a simply connected symplectic
4-manifold $W$ with
$b^+=1$, $b^-=2$, and a symplectically embedded sphere representing
$(2p-7)\a+\b$,
a class of square $2p-7>0$. As above, McDuff's result implies that
$W=\CP\#2\CPb$.
The class $\z$ is represented by an exceptional sphere, which we now
blow down to
obtain a manifold $Y$, which must be diffeomorphic to either
$\CP\#\CPb$ or $\SS$.
The complement of the rational ball $B_{2p-7}$ in $Y$ has its second homology
generated by the classes $(2p-7)\a+\b+\z$ and $\a-\ve$ with intersection 
matrix:
\[ \begin{pmatrix} 2p-6 & 1 \\ 1 & -(2p-8)
\end{pmatrix}\] Thus $Y=\Bbb{F}_{2p-8}\cong\SS$ and $s_++f=(2p-7)\a+\b+\z$,
$s_-=\a-\ve$.

Let $A, B$ denote the classes $[S^2\x \{\text{pt}\}]$,
$[\{\text{pt}\}\x S^2]$ in
$H_2(\SS)$ where the fiber $f=B$. Then if we identify
$Y\#\CPb=W=\CP\#2\CPb$ this
identifies:
$A\lr h-e_1$, $B\lr h-e$, and $\z\lr h-e-e_1$.
Now $\Sig'_R$ has been blown down in $W$ to represent
\begin{multline*} \z+(\a-\ve)+2((2p-7)\a+\b) \ = \  2(s_++f)+s_--\z\\
=2(A+(p-3)B)+(A-(p-4)B)-\z \ = \ p\,h-(p-3)\,e-2\,e_1
\end{multline*}
which is how $\Sig'_S$ is constructed.
\end{proof}

The rational surface $R(q+1)=\CP\#\CPb\#\,4q\,\CPb$ may be obtained as the
(desingularized) double branched cover of $\SS$, branched over two copies of
$S^2\x\{\text{pt}\}$ and $2q$ copies of $\{\text{pt}\}\x S^2$. In
this way we see
that $R(q+1)$ admits a `vertical' genus $0$  fibration over $S^2$
with fiber class
$H-E$ and also a genus $q-1$ `horizontal' fibration over $S^2$.

\begin{lem} The fibers of the horizontal fibration of $R(q+1)$ are isotopic to
$\Sig_{R(q+1)}$.\end{lem}
\begin{proof} The vertical fibration on $R(q+1)$ has $2q$ singular fibers, each
consisting of an exceptional sphere of multiplicity $2$
together with a pair of disjoint spheres of self-intersection $-2$, each
intersecting the exceptional sphere in a single point. Consider the first such
singular fiber --- call the spheres, $E_1$, $x$, and $y$.  Blowing down $E_1$
leaves a pair of exceptional curves, $x+E_1$ and $y+E_1$. Blow down
$E_2 = x+E_1$ to
obtain a single sphere $y+E_1+E_2$ whose square is $0$. This is now the fiber
$H-E$ of a genus $0$ fibration of $\CP\#\CPb\#\,(4q-2)\,\CPb$. It follows that
$y=H-E-E_1-E_2$. In general, the $i$th singular fiber of the vertical
fibration on
$R(q+1)$ consists of an exceptional curve $E_{2i-1}$ of multiplicity
$2$, along with
a pair of disjoint ($-2$)-spheres, $E_{2i}-E_{2i-1}$ and $H-E-E_{2i-1}-E_{2i}$.

The horizontal fiber $\L$ is homologous to $aH-bE-\sum_1^{4q}c_iE_i$ for some
coefficients $a$, $b$, $c_i$. Since a generic horizontal fiber
intersects a generic
vertical fiber in two points,
$2=\L\cdot (H-E)= a-b$. Furthermore, a generic horizontal fiber is
disjoint from
the ($-2$)-spheres which occur as part of the vertical singular fibers. Thus
$\L\cdot (E_{2i}-E_{2i-1}) = 0 = \L\cdot (H-E-E_{2i-1}-E_{2i})$. The
first of these
two equalities shows that $c_{2i-1}=c_{2i}$ for $i=1,\dots,2q$. The
second shows
that $a-(a-2)-c_{2i-1}-c_{2i}=0$; so $c_{2i-1}+c_{2i}=2$, and thus $c_i=1$,
$i=1,\dots,4q$. Finally, $\L^2=0$ gives $a=q+1$. This shows that $\L$ and
$\Sig_{R(q+1)}$ are homologous. Since both are holomorphic curves in
$R(q+1)$, they
must actually be isotopic (cf.\ the introduction to \cite{surfaces}).
\end{proof}

We can now calculate the Seiberg-Witten invariants of $X_p$ and
$X'_p$. Let $E(q)$
denote the simply connected elliptic surface with $\ch=q$ and with no multiple
fibers.

\begin{lem} $E(q)$ is diffeomorphic to the fiber sum $R(q+1)\#_{\Sig_{R(q+1)}}
R(q+1)$. \end{lem}
\begin{proof} The (desingularized) double cover of $\SS$ branched
over four copies of
$S^2\x\{\text{pt}\}$ and $2q$ copies of $\{\text{pt}\}\x S^2$ is $E(q)$. The
previous lemma shows that this is the fiber sum as advertised.
\end{proof}

It follows from this lemma and Proposition~\ref{P} that $X_p$ is the rational
blowdown of a configuration $C_{2p-6}$ in $E(2p-4)$. The elliptic fiber $T$
of \[ E(2p-4)=R(2p-3)\#_{\Sig_{R(2p-3)}}R(2p-3)\] is obtained from a
genus zero fiber
on each side, since these spheres intersect $\Sig_{R(2p-3)}$ in two points. The
genus 0 fiber in $R(2p-3)$ represents the class $H-E$, and the lead
sphere $S_0$ of
$C_{2p-6}$ represents $H-\sum_1^{2p-3}E_i$.

The basic classes of $E(2p-4)$ are $\pm 2j\,T$, $j=0,\dots,p-3$. Of these, only
$\pm(2p-6)\,T$ intersects $S_0$ maximally (with intersection number
$\pm(2p-6)$). It
follows from Theorem~\ref{taut} that the rational blowdown $X_p$ has
(up to sign)
just one basic class. A similar argument shows that the same is true
for $X'_p$.

\begin{prop} The simply connected symplectic manifolds $X_p$ and
$X'_p$\break ($p\ge4$)
have (up to sign) one basic class and satisfy $c_1^2=\ch -3$.
\qed\end{prop}

\section{Construction 1}

In order to fill in the region, $\ch - 3\le c_1^2 \le 2\ch - 6$
we shall next exhibit symplectic spheres of self-intersection $-4$ in
$X(p)$ and
$X'(p)$ which can be rationally blown down. These spheres will be built from
`pieces' which intersect $\Sig_{R(2p-3)}$ and $\Sig_{S(p)}$ transversely. In
$R(2p-3)$ there are the exceptional spheres $E_i$, $i=1,2,\dots, 8p-16$ which
intersect $\Sig_{R(2p-3)}$ transversely in a single point. Also,
consider a line
in $\CP$ that does not belong to the arrangement for $R(2p-3)$, but which goes
through the singular point of order $2p-5$ in the arrangement.
This line gets blown up to a sphere of self-intersection
$0$ which intersects $\Sig_{R(2p-3)}$ transversely in two points. This sphere
represents the class $H-E$. We may form arbitrarily many such disjoint spheres.
Denote them by $A_j$.

To construct spheres in $S(p)$, recall how it is constructed. The initial
arrangement consists of $p-3$ lines through a common point $x_0$, and
three further
lines $L_i$, $i=1,\dots,3$, in general position. One then blows up at
$x_0$ and at
further points $x_j$, $j=1,\dots,6p-9$ on the arrangement. We can suppose the
$x_{2k-1}$  lie on $L_1$ and that $x_{2k}$  lie on $L_2$ for
$k=1,\dots,3p-5$ and are arranged so that each pair of points
$\{x_{2k-1},x_{2k}\}$
lies on a line $B_k'$ through $x_0$. After all the blowups, one obtains spheres
$B_k$ ($k=1,\dots,3p-5$) of self-intersection $-2$ in $S(p)$. These
spheres intersect
$\Sig_{S(p)}$ transversely in one point (the point of intersection of
$B_k'$ with
$L_3$). Note that $B_k$ is homologous to $H-E-E_{2k-1}-E_{2k}$.

Also, there
are spheres $C_{\ell}$ with self-intersection $0$ that intersect
$\Sig_{S(p)}$ in
three points that are obtained from a line in $\CP$ that goes through
the singular
point of order $p-3$. The spheres $C_{\ell}$ are homologous to $H-E$.

In $X_p$ each of the spheres $E_i$, $A_j$, $B_k$, and $C_{\ell}$ is punctured.
One can form the fiber sum so that the punctures match up in such a
way that $B_1
\cup A_1 \cup C_1\cup E_1 \cup E_2$ is a symplectic sphere of
self-intersection $-4$
in $X_p$. Further, we can arrange so that there are $3p-5$ disjoint
symplectic ($-4$)-spheres constructed in this way. Rationally blowing
these down, one
at a time, we obtain simply connected symplectic manifolds $X(p,k)$
which have, up to
sign one basic class, and with $\ch(X(p,k))=\ch(X_p)= 2p-4$, and
$c_1^2(X(p,k))=c_1^2(X_p)+k= 2p-7+k$, i.e. filling up the region
$\ch-3\le c_1^2\le
\frac{5}{2}\ch -2$, $\ch$ even.

The same construction applied to $X'_p$ yields the odd $\ch$
examples. In this case
one can construct $3p-7$ of the spheres $B_k$ and thence $3p-7$ spheres of
self-intersection $-4$ to rationally blow down. We get manifolds $X'(p,k)$ with
$\ch(X'(p,k))=\ch(X_p)= 2p-5$, and $c_1^2(X(p,k))=c_1^2(X'_p)+k=
2p-8+k$. So we fill the region $\ch-3\le c_1^2\le \frac{5}{2}\ch -2$, $\ch$
odd.

\begin{thm}\label{TT} For every pair of positive integers $(x,c)$
with $0< x-3 \le c
\le \frac{5}{2}x-2$ there is a simply connected symplectic 4-manifold $X$ with
$c_1^2(X)=c$,
$\chi(X)=x$ and (up to sign) one basic class.
\end{thm}

This implies Theorem~\ref{T}.

\section{Construction 2}

We shall now give a second proof of Theorem~\ref{TT} with a
construction starting
directly with the elliptic surfaces $E(n)$. Fix a pair of positive
integers $(x,c)$
with $0 < x-3 \le c \le \frac{5}{2}x-2$ as in the statement of the theorem, and
consider the elliptic surface $E(x)$. It admits an elliptic fibration with $6x$
cusp fibers and no other singular fibers. Furthermore, $E(x)$ contains, as a
symplectic codimension $0$ submanifold, the canonical resolution of the $(2,
2x-1, 4x-3)$ Brieskorn singularity. This contains the configuration
of symplectic
spheres:

\vspace{-3mm}\centerline{\unitlength .9cm\small
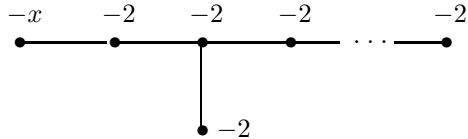
\begin{picture}(6,2)
\put(-.5,.7){$\bullet$}
\put(-.4,.8){\line(1,0){1.3}}
\put(.9,.7){$\bullet$}
\put(1,.8){\line(1,0){1.3}}
\put(2.2,.7){$\bullet$}
\put(2.3,.8){\line(1,0){1.3}}
\put(3.5,.7){$\bullet$}
\put(3.6,.8){\line(1,0){.75}}
\put(4.55,.8){.}
\put(4.75,.8){.}
\put(4.95,.8){.}
\put(5.15,.8){\line(1,0){.75}}
\put(5.8,.7){$\bullet$}
\put(2.3,.8){\line(0,-1){1.3}}
\put(2.2,-.6){$\bullet$}
\put(-.6,1.1){$-x$}
\put(.8,1.1){$-2$}
\put(2.1,1.1){$-2$}
\put(3.4,1.1){$-2$}
\put(5.7,1.1){$-2$}
\put(2.5,-.6){$-2$}
\end{picture}}

\vspace{.2in}
\noindent where the linear plumbing to the right of the central node has $4x-4$
spheres of self-intersection $-2$, and where the sphere $S$ of
self-intersection
$-x$ is a section of the elliptic fibration on $E(x)$  (cf.\ \cite{canonical}).

Each fiber of the elliptic fibration meets $S$ in a single positive
intersection.
In particular, consider one of the $6x$ cusp fibers. If we blow up at
the cuspidal
point, then in $E(x)\# \CPb$ we obtain an embedded sphere
representing $f-2E$ where
$f$ is the class of an elliptic fiber in $E(x)$ and $E$ is the class of the
exceptional curve. Thus $f-2E$ represents an embedded sphere that
meets $S$ once,
positively, in $E(x)\# \CPb$. Now symplectically resolve the
intersection to obtain
a sphere $S'$ of self-intersection $-(x+2)$. Comparing with the
plumbing diagram
above, we see a symplectic embedding of the configuration $C_x$ in
$E(x)\# \CPb$.
This process can be repeated until we either exhaust all $6x$ of the
cusp fibers
or all $4x-2$ of the  ($-2$)-spheres across the top of the plumbing.

If we blow up $k$ of the cuspidal points, we obtain a sphere of
self-intersection $-(x+2k)$. This is the lead sphere of the 
configuration $C_{x+2k-2}$
which has $x+2k-4$ ($-2$)-spheres. Thus we can find a symplectic configuration
$C_{x+2k-2}$ in $E(x)\#\,k\,\CPb$ for $x+2k-4\le 4x-2$; i.e. for $0 \le k \le
\frac{3}{2} x+1$. (The $k=0$ case is $C_{x-2}\C E(x)$.)

In $E(x)\#\,k\,\CPb$, rationally blow down the configuration $C_{x+2k-2}$. This
yields a symplectic manifold with $\ch = x$ and $c_1^2 = x+k-3$.
Applying the blowup
formula \cite{Turkey} shows that the Seiberg-Witten basic
classes of $E(x)\#\,k\,\CPb$ have the form
$$\b(m;\ve_1,\dots,\ve_k) = mf+\ve_1E_1+\dots +\ve_kE_k$$ for $|m|\le
x-2$, $m \equiv
x\pmod2$, and $\ve_i=\pm1$ for $1\le i\le k$. We now apply
Theorem~\ref{taut}. Each
of the ($-2$)-spheres $S_i$, $i=1,\dots, x+2k-4$, is embedded in $E(x)$ and is
symplectic. Therefore, it follows from the adjunction formula that
$$ 0= K_{E(x)}\cdot S_i = (x-2)f\cdot S_i $$
where $K_{E(x)}$ denotes the canonical class of $E(x)$. Hence $f\cdot
S_i=0$ for
$i=1,\dots, x+2k-4$. Since clearly each
$E_j\cdot S_i=0$, we have
\[ \b(m;\ve_1,\dots,\ve_k)\cdot S_i = 0, \ \ i=1,\dots, x+2k-4.\]
The lead sphere $S_0$ of our configuration $C_{x-2k-2}$ is given
homologically by
$$S_0=S-2E_1-\dots-2E_k. $$
$$\leqno{\rm Hence}\qquad\b(m;\ve_1,\dots,\ve_k)\cdot S_0 = m + 2\sum_{i=1}^k\ve_i \le m+2k
\le x+2k-2.$$

Thus the hypotheses of Theorem~\ref{taut} are satisfied.
It is now easy to see that only the basic classes $\pm\b(x-2;1,\dots,1)$
satisfy $\b(m;\ve_1,\dots,\ve_k)\cdot S_0 = x+2k-2$; and so, up to sign, our
manifold has just one basic class.

Since this construction yields 4-manifolds with $\ch=x$ and
$c_1^2 = x+k-3$ for $0 \le k\le \frac{3}{2}x +1$, the existence of
these manifolds
again proves Theorem~\ref{TT}.

The authors do not know if the families of manifolds produced by our two
constructions actually coincide. This is quite plausible and seems to be an
interesting question.


\Addresses\recd

\end{document}

%% file: agtout.tex
\def\ifplaintex{\expandafter\ifx\csname documentclass\endcsname\relax}

\def\gtp{{\mathsurround=0pt\it $\cal G\mskip-2mu$eometry \&\ 
$\cal T\!\!$opology $\cal P\!$ublications}}  

\def\recd{{\small Received:\qua\receiveddate\ifx\reviseddate\relax
\else\qquad Revised:\qua\reviseddate\fi\par}} 

\def\lognumber#1{\def\thelognumber{#1}}
\def\volumenumber#1{\def\thevolumenumber{#1}}
\def\volumeyear#1{\def\thevolumeyear{#1}}
\def\papernumber#1{\def\thepapernumber{#1}}
\def\pagenumbers#1#2{\def\startpage{#1}\def\finishpage{#2}}
\def\published#1{\def\publishdate{#1}}

\def\received#1{\def\receiveddate{#1}}

\def\accepted#1{\def\accepteddate{#1}}

\def\asciiaddress#1{\def\theasciiaddress{#1}}

\let\\\par\let\thelognumber\relax\let\thevolumenumber\relax
\let\thepapernumber\relax\let\thevolumeyear\relax\let\startpage\relax
\let\finishpage\relax\let\publishdate\relax\let\receiveddate\relax
\let\reviseddate\relax\let\accepteddate\relax\let\theasciititle\relax
\let\theasciiauthors\relax\let\theasciiaddress\relax
\let\theasciiabstract\relax

\let\theasciiemail\relax

\ifplaintex
\font\logobig=cmssbx10 scaled 3836
\font\logomed=cmssbx10 scaled 2557
\else
\font\logobig=cmssbx10 scaled 4200
\font\logomed=cmssbx10 scaled 2800
\fi

\long\def\makeagttitle{   
\count0=\startpage
\agt\hfill      
\hbox to 45truept{\vbox to 0pt{\vglue -13truept{\logomed A\kern -.37em{\logobig 
T}\kern -.38em G}\vss}\hss}
\break
{\small Volume \thevolumenumber\ (\thevolumeyear)
\startpage--\finishpage\nl
Published: \publishdate}

\vglue .25truein

{\parskip=0pt\leftskip 0pt plus
1fil\def\\{\par\smallskip}{\Large\bf\thetitle}\par\medskip} \vglue
0.05truein

{\parskip=0pt\leftskip 0pt plus 1fil\def\\{\par}{\sc\theauthors}
\par\medskip}%
 
\vglue 0.03truein 

{\small\leftskip 25truept\rightskip 25truept{\bf Abstract}\stdspace\theabstract

{\bf AMS Classification}\stdspace\theprimaryclass
\ifx\thesecondaryclass\relax\else; \thesecondaryclass\fi\par
{\bf Keywords}\stdspace \thekeywords\par}\vglue 7truept

}   

\ifplaintex
\font\phead=cmsl9 scaled 950
\font\pnum=cmbx10 scaled 913
\font\pfoot=cmsl9 scaled 950
\headline{\vbox to 0pt{\vskip -4.5mm\line{\small\phead\ifnum
\count0=\startpage ISSN 1472-2739 (on-line) 1472-2747 (printed)
\hfill {\pnum\folio}\else\ifodd\count0\def\\{ }%
\ifx\theshorttitle\relax\thetitle\else\theshorttitle\fi\hfill{\pnum\folio}
\else\def\\{ and }{\pnum\folio}\hfill\ifx\theshortauthors\relax\theauthors
\else\theshortauthors\fi\fi\fi}\vss}}
\footline{\vbox to 0pt{\vglue 0mm\line{\small\pfoot\ifnum\count0=\startpage
\copyright\ \gtp\hfill\else
\agt, Volume \thevolumenumber\ (\thevolumeyear)\hfill\fi}\vss}}
\else
\font=cmsl9 scaled 1050
\font\lnum=cmbx10 
\font=cmsl9 scaled 1050
\makeatletter
\def\@oddhead{{\small\lhead\ifnum\count0=\startpage ISSN 1472-2739 
(on-line) 1472-2747 (printed)\hfill {\lnum\number\count0}\else\ifodd\count0
\def\\{ }\ifx\theshorttitle\relax \thetitle \else\theshorttitle\fi\hfill
{\lnum\number\count0}\else\def\\{ and }{\lnum\number\count0}
\hfill\ifx\theshortauthors\relax 
\theauthors\else\theshortauthors\fi\fi\fi}}\def\@evenhead{\@oddhead}
\def\@oddfoot{\small\lfoot\ifnum\count0=\startpage\copyright\ \gtp\hfill\else
\agt, Volume \thevolumenumber\ (\thevolumeyear)\hfill\fi}
\def\@evenfoot{\@oddfoot}
\makeatother
\fi
\let\maketitlepage\makeagttitle
\let\makeshorttitle\maketitlepage
\let\maketitle\maketitlepage

\newwrite\gtoutfile
\long\gdef\makeheadfile{  
{\def\\{, }\def\s{ }
\immediate\openout\gtoutfile head.xxx
\immediate\write\gtoutfile{To: math@arxiv.org}
\immediate\write\gtoutfile{Subject: put OR rep NNNNN:ppppp}
\immediate\write\gtoutfile{--text follows this line--}
\immediate\write\gtoutfile{Proxy-for: \ifx\theasciiauthors\relax
\theauthors\else\theasciiauthors\fi\s<\ifx\theasciiemail\relax\theemail\else\theasciiemail\fi>}
\immediate\write\gtoutfile{\noexpand\\}
\immediate\write\gtoutfile{Authors: \ifx\theasciiauthors\relax
\theauthors\else\theasciiauthors\fi}
{\def\\{ }\immediate\write\gtoutfile{Title: \ifx\theasciititle\relax
\thetitle\else\theasciititle\fi}}
\immediate\write\gtoutfile{Subj-class: GT or SG, GR etc}
\immediate\write\gtoutfile{MSC-class: \theprimaryclass\ifx\thesecondaryclass\relax\else, \thesecondaryclass\fi}
\immediate\write\gtoutfile{Journal-ref: Algebr. Geom. Topol. \thevolumenumber\s
(\thevolumeyear) \startpage-\finishpage}
\immediate\write\gtoutfile{Comments: Published by Algebraic and
Geometric Topology at}
\immediate\write\gtoutfile{\s\s\s  http://www.maths.warwick.ac.uk/agt/AGTVol\thevolumenumber/agt-\thevolumenumber-\thepapernumber.abs.html}
\immediate\write\gtoutfile{\noexpand\\}
\immediate\write\gtoutfile{}
\ifx\theasciiabstract\relax
\immediate\write\gtoutfile{\theabstract}\else
\immediate\write\gtoutfile{\theasciiabstract}\fi
\immediate\write\gtoutfile{}
\immediate\write\gtoutfile{\noexpand\\}
\immediate\write\gtoutfile{}
\immediate\closeout\gtoutfile}}  

\def\maketitlepage{\makeagttitle\makeheadfile}
\let\makeshorttitle\maketitlepage
\let\maketitle\maketitlepage

\def\ifplaintex{\expandafter\ifx\csname documentclass\endcsname\relax}

\def\gtp{{\mathsurround=0pt\it $\cal G\mskip-2mu$eometry \&\ 
$\cal T\!\!$opology $\cal P\!$ublications}}  

\def\recd{{\small Received:\qua\receiveddate\ifx\reviseddate\relax
\else\qquad Revised:\qua\reviseddate\fi\par}} 

\def\lognumber#1{\def\thelognumber{#1}}
\def\volumenumber#1{\def\thevolumenumber{#1}}
\def\volumeyear#1{\def\thevolumeyear{#1}}
\def\papernumber#1{\def\thepapernumber{#1}}
\def\pagenumbers#1#2{\def\startpage{#1}\def\finishpage{#2}}
\def\published#1{\def\publishdate{#1}}

\def\received#1{\def\receiveddate{#1}}

\def\accepted#1{\def\accepteddate{#1}}

\def\asciiaddress#1{\def\theasciiaddress{#1}}

\let\\\par\let\thelognumber\relax\let\thevolumenumber\relax
\let\thepapernumber\relax\let\thevolumeyear\relax\let\startpage\relax
\let\finishpage\relax\let\publishdate\relax\let\receiveddate\relax
\let\reviseddate\relax\let\accepteddate\relax\let\theasciititle\relax
\let\theasciiauthors\relax\let\theasciiaddress\relax
\let\theasciiabstract\relax

\let\theasciiemail\relax

\ifplaintex
\font\logobig=cmssbx10 scaled 3836
\font\logomed=cmssbx10 scaled 2557
\else
\font\logobig=cmssbx10 scaled 4200
\font\logomed=cmssbx10 scaled 2800
\fi

\long\def\makeagttitle{   
\count0=\startpage
\agt\hfill      
\hbox to 45truept{\vbox to 0pt{\vglue -13truept{\logomed A\kern -.37em{\logobig 
T}\kern -.38em G}\vss}\hss}
\break
{\small Volume \thevolumenumber\ (\thevolumeyear)
\startpage--\finishpage\nl
Published: \publishdate}

\vglue .25truein

{\parskip=0pt\leftskip 0pt plus
1fil\def\\{\par\smallskip}{\Large\bf\thetitle}\par\medskip} \vglue
0.05truein

{\parskip=0pt\leftskip 0pt plus 1fil\def\\{\par}{\sc\theauthors}
\par\medskip}%
 
\vglue 0.03truein 

{\small\leftskip 25truept\rightskip 25truept{\bf Abstract}\stdspace\theabstract

{\bf AMS Classification}\stdspace\theprimaryclass
\ifx\thesecondaryclass\relax\else; \thesecondaryclass\fi\par
{\bf Keywords}\stdspace \thekeywords\par}\vglue 7truept

}   

\ifplaintex
\font\phead=cmsl9 scaled 950
\font\pnum=cmbx10 scaled 913
\font\pfoot=cmsl9 scaled 950
\headline{\vbox to 0pt{\vskip -4.5mm\line{\small\phead\ifnum
\count0=\startpage ISSN 1472-2739 (on-line) 1472-2747 (printed)
\hfill {\pnum\folio}\else\ifodd\count0\def\\{ }%
\ifx\theshorttitle\relax\thetitle\else\theshorttitle\fi\hfill{\pnum\folio}
\else\def\\{ and }{\pnum\folio}\hfill\ifx\theshortauthors\relax\theauthors
\else\theshortauthors\fi\fi\fi}\vss}}
\footline{\vbox to 0pt{\vglue 0mm\line{\small\pfoot\ifnum\count0=\startpage
\copyright\ \gtp\hfill\else
\agt, Volume \thevolumenumber\ (\thevolumeyear)\hfill\fi}\vss}}
\else
\font=cmsl9 scaled 1050
\font\lnum=cmbx10 
\font=cmsl9 scaled 1050
\makeatletter
\def\@oddhead{{\small\lhead\ifnum\count0=\startpage ISSN 1472-2739 
(on-line) 1472-2747 (printed)\hfill {\lnum\number\count0}\else\ifodd\count0
\def\\{ }\ifx\theshorttitle\relax \thetitle \else\theshorttitle\fi\hfill
{\lnum\number\count0}\else\def\\{ and }{\lnum\number\count0}
\hfill\ifx\theshortauthors\relax 
\theauthors\else\theshortauthors\fi\fi\fi}}\def\@evenhead{\@oddhead}
\def\@oddfoot{\small\lfoot\ifnum\count0=\startpage\copyright\ \gtp\hfill\else
\agt, Volume \thevolumenumber\ (\thevolumeyear)\hfill\fi}
\def\@evenfoot{\@oddfoot}
\makeatother
\fi
\let\maketitlepage\makeagttitle
\let\makeshorttitle\maketitlepage
\let\maketitle\maketitlepage

\newwrite\gtoutfile
\long\gdef\makeheadfile{  
{\def\\{, }\def\s{ }
\immediate\openout\gtoutfile head.xxx
\immediate\write\gtoutfile{To: math@arxiv.org}
\immediate\write\gtoutfile{Subject: put OR rep NNNNN:ppppp}
\immediate\write\gtoutfile{--text follows this line--}
\immediate\write\gtoutfile{Proxy-for: \ifx\theasciiauthors\relax
\theauthors\else\theasciiauthors\fi\s<\ifx\theasciiemail\relax\theemail\else\theasciiemail\fi>}
\immediate\write\gtoutfile{\noexpand\\}
\immediate\write\gtoutfile{Authors: \ifx\theasciiauthors\relax
\theauthors\else\theasciiauthors\fi}
{\def\\{ }\immediate\write\gtoutfile{Title: \ifx\theasciititle\relax
\thetitle\else\theasciititle\fi}}
\immediate\write\gtoutfile{Subj-class: GT or SG, GR etc}
\immediate\write\gtoutfile{MSC-class: \theprimaryclass\ifx\thesecondaryclass\relax\else, \thesecondaryclass\fi}
\immediate\write\gtoutfile{Journal-ref: Algebr. Geom. Topol. \thevolumenumber\s
(\thevolumeyear) \startpage-\finishpage}
\immediate\write\gtoutfile{Comments: Published by Algebraic and
Geometric Topology at}
\immediate\write\gtoutfile{\s\s\s  http://www.maths.warwick.ac.uk/agt/AGTVol\thevolumenumber/agt-\thevolumenumber-\thepapernumber.abs.html}
\immediate\write\gtoutfile{\noexpand\\}
\immediate\write\gtoutfile{}
\ifx\theasciiabstract\relax
\immediate\write\gtoutfile{\theabstract}\else
\immediate\write\gtoutfile{\theasciiabstract}\fi
\immediate\write\gtoutfile{}
\immediate\write\gtoutfile{\noexpand\\}
\immediate\write\gtoutfile{}
\immediate\closeout\gtoutfile}}  

\def\maketitlepage{\makeagttitle\makeheadfile}
\let\makeshorttitle\maketitlepage
\let\maketitle\maketitlepage

\def\ifplaintex{\expandafter\ifx\csname documentclass\endcsname\relax}

\def\gtp{{\mathsurround=0pt\it $\cal G\mskip-2mu$eometry \&\ 
$\cal T\!\!$opology $\cal P\!$ublications}}  

\def\recd{{\small Received:\qua\receiveddate\ifx\reviseddate\relax
\else\qquad Revised:\qua\reviseddate\fi\par}} 

\def\lognumber#1{\def\thelognumber{#1}}
\def\volumenumber#1{\def\thevolumenumber{#1}}
\def\volumeyear#1{\def\thevolumeyear{#1}}
\def\papernumber#1{\def\thepapernumber{#1}}
\def\pagenumbers#1#2{\def\startpage{#1}\def\finishpage{#2}}
\def\published#1{\def\publishdate{#1}}

\def\received#1{\def\receiveddate{#1}}

\def\accepted#1{\def\accepteddate{#1}}

\def\asciiaddress#1{\def\theasciiaddress{#1}}

\let\\\par\let\thelognumber\relax\let\thevolumenumber\relax
\let\thepapernumber\relax\let\thevolumeyear\relax\let\startpage\relax
\let\finishpage\relax\let\publishdate\relax\let\receiveddate\relax
\let\reviseddate\relax\let\accepteddate\relax\let\theasciititle\relax
\let\theasciiauthors\relax\let\theasciiaddress\relax
\let\theasciiabstract\relax

\let\theasciiemail\relax

\ifplaintex
\font\logobig=cmssbx10 scaled 3836
\font\logomed=cmssbx10 scaled 2557
\else
\font\logobig=cmssbx10 scaled 4200
\font\logomed=cmssbx10 scaled 2800
\fi

\long\def\makeagttitle{   
\count0=\startpage
\agt\hfill      
\hbox to 45truept{\vbox to 0pt{\vglue -13truept{\logomed A\kern -.37em{\logobig 
T}\kern -.38em G}\vss}\hss}
\break
{\small Volume \thevolumenumber\ (\thevolumeyear)
\startpage--\finishpage\nl
Published: \publishdate}

\vglue .25truein

{\parskip=0pt\leftskip 0pt plus
1fil\def\\{\par\smallskip}{\Large\bf\thetitle}\par\medskip} \vglue
0.05truein

{\parskip=0pt\leftskip 0pt plus 1fil\def\\{\par}{\sc\theauthors}
\par\medskip}%
 
\vglue 0.03truein 

{\small\leftskip 25truept\rightskip 25truept{\bf Abstract}\stdspace\theabstract

{\bf AMS Classification}\stdspace\theprimaryclass
\ifx\thesecondaryclass\relax\else; \thesecondaryclass\fi\par
{\bf Keywords}\stdspace \thekeywords\par}\vglue 7truept

}   

\ifplaintex
\font\phead=cmsl9 scaled 950
\font\pnum=cmbx10 scaled 913
\font\pfoot=cmsl9 scaled 950
\headline{\vbox to 0pt{\vskip -4.5mm\line{\small\phead\ifnum
\count0=\startpage ISSN 1472-2739 (on-line) 1472-2747 (printed)
\hfill {\pnum\folio}\else\ifodd\count0\def\\{ }%
\ifx\theshorttitle\relax\thetitle\else\theshorttitle\fi\hfill{\pnum\folio}
\else\def\\{ and }{\pnum\folio}\hfill\ifx\theshortauthors\relax\theauthors
\else\theshortauthors\fi\fi\fi}\vss}}
\footline{\vbox to 0pt{\vglue 0mm\line{\small\pfoot\ifnum\count0=\startpage
\copyright\ \gtp\hfill\else
\agt, Volume \thevolumenumber\ (\thevolumeyear)\hfill\fi}\vss}}
\else
\font=cmsl9 scaled 1050
\font\lnum=cmbx10 
\font=cmsl9 scaled 1050
\makeatletter
\def\@oddhead{{\small\lhead\ifnum\count0=\startpage ISSN 1472-2739 
(on-line) 1472-2747 (printed)\hfill {\lnum\number\count0}\else\ifodd\count0
\def\\{ }\ifx\theshorttitle\relax \thetitle \else\theshorttitle\fi\hfill
{\lnum\number\count0}\else\def\\{ and }{\lnum\number\count0}
\hfill\ifx\theshortauthors\relax 
\theauthors\else\theshortauthors\fi\fi\fi}}\def\@evenhead{\@oddhead}
\def\@oddfoot{\small\lfoot\ifnum\count0=\startpage\copyright\ \gtp\hfill\else
\agt, Volume \thevolumenumber\ (\thevolumeyear)\hfill\fi}
\def\@evenfoot{\@oddfoot}
\makeatother
\fi
\let\maketitlepage\makeagttitle
\let\makeshorttitle\maketitlepage
\let\maketitle\maketitlepage

\newwrite\gtoutfile
\long\gdef\makeheadfile{  
{\def\\{, }\def\s{ }
\immediate\openout\gtoutfile head.xxx
\immediate\write\gtoutfile{To: math@arxiv.org}
\immediate\write\gtoutfile{Subject: put OR rep NNNNN:ppppp}
\immediate\write\gtoutfile{--text follows this line--}
\immediate\write\gtoutfile{Proxy-for: \ifx\theasciiauthors\relax
\theauthors\else\theasciiauthors\fi\s<\ifx\theasciiemail\relax\theemail\else\theasciiemail\fi>}
\immediate\write\gtoutfile{\noexpand\\}
\immediate\write\gtoutfile{Authors: \ifx\theasciiauthors\relax
\theauthors\else\theasciiauthors\fi}
{\def\\{ }\immediate\write\gtoutfile{Title: \ifx\theasciititle\relax
\thetitle\else\theasciititle\fi}}
\immediate\write\gtoutfile{Subj-class: GT or SG, GR etc}
\immediate\write\gtoutfile{MSC-class: \theprimaryclass\ifx\thesecondaryclass\relax\else, \thesecondaryclass\fi}
\immediate\write\gtoutfile{Journal-ref: Algebr. Geom. Topol. \thevolumenumber\s
(\thevolumeyear) \startpage-\finishpage}
\immediate\write\gtoutfile{Comments: Published by Algebraic and
Geometric Topology at}
\immediate\write\gtoutfile{\s\s\s  http://www.maths.warwick.ac.uk/agt/AGTVol\thevolumenumber/agt-\thevolumenumber-\thepapernumber.abs.html}
\immediate\write\gtoutfile{\noexpand\\}
\immediate\write\gtoutfile{}
\ifx\theasciiabstract\relax
\immediate\write\gtoutfile{\theabstract}\else
\immediate\write\gtoutfile{\theasciiabstract}\fi
\immediate\write\gtoutfile{}
\immediate\write\gtoutfile{\noexpand\\}
\immediate\write\gtoutfile{}
\immediate\closeout\gtoutfile}}  

\def\maketitlepage{\makeagttitle\makeheadfile}
\let\makeshorttitle\maketitlepage
\let\maketitle\maketitlepage

\def\ifplaintex{\expandafter\ifx\csname documentclass\endcsname\relax}

\def\gtp{{\mathsurround=0pt\it $\cal G\mskip-2mu$eometry \&\ 
$\cal T\!\!$opology $\cal P\!$ublications}}  

\def\recd{{\small Received:\qua\receiveddate\ifx\reviseddate\relax
\else\qquad Revised:\qua\reviseddate\fi\par}} 

\def\lognumber#1{\def\thelognumber{#1}}
\def\volumenumber#1{\def\thevolumenumber{#1}}
\def\volumeyear#1{\def\thevolumeyear{#1}}
\def\papernumber#1{\def\thepapernumber{#1}}
\def\pagenumbers#1#2{\def\startpage{#1}\def\finishpage{#2}}
\def\published#1{\def\publishdate{#1}}

\def\received#1{\def\receiveddate{#1}}

\def\accepted#1{\def\accepteddate{#1}}

\def\asciiaddress#1{\def\theasciiaddress{#1}}

\let\\\par\let\thelognumber\relax\let\thevolumenumber\relax
\let\thepapernumber\relax\let\thevolumeyear\relax\let\startpage\relax
\let\finishpage\relax\let\publishdate\relax\let\receiveddate\relax
\let\reviseddate\relax\let\accepteddate\relax\let\theasciititle\relax
\let\theasciiauthors\relax\let\theasciiaddress\relax
\let\theasciiabstract\relax

\let\theasciiemail\relax

\ifplaintex
\font\logobig=cmssbx10 scaled 3836
\font\logomed=cmssbx10 scaled 2557
\else
\font\logobig=cmssbx10 scaled 4200
\font\logomed=cmssbx10 scaled 2800
\fi

\long\def\makeagttitle{   
\count0=\startpage
\agt\hfill      
\hbox to 45truept{\vbox to 0pt{\vglue -13truept{\logomed A\kern -.37em{\logobig 
T}\kern -.38em G}\vss}\hss}
\break
{\small Volume \thevolumenumber\ (\thevolumeyear)
\startpage--\finishpage\nl
Published: \publishdate}

\vglue .25truein

{\parskip=0pt\leftskip 0pt plus
1fil\def\\{\par\smallskip}{\Large\bf\thetitle}\par\medskip} \vglue
0.05truein

{\parskip=0pt\leftskip 0pt plus 1fil\def\\{\par}{\sc\theauthors}
\par\medskip}%
 
\vglue 0.03truein 

{\small\leftskip 25truept\rightskip 25truept{\bf Abstract}\stdspace\theabstract

{\bf AMS Classification}\stdspace\theprimaryclass
\ifx\thesecondaryclass\relax\else; \thesecondaryclass\fi\par
{\bf Keywords}\stdspace \thekeywords\par}\vglue 7truept

}   

\ifplaintex
\font\phead=cmsl9 scaled 950
\font\pnum=cmbx10 scaled 913
\font\pfoot=cmsl9 scaled 950
\headline{\vbox to 0pt{\vskip -4.5mm\line{\small\phead\ifnum
\count0=\startpage ISSN 1472-2739 (on-line) 1472-2747 (printed)
\hfill {\pnum\folio}\else\ifodd\count0\def\\{ }%
\ifx\theshorttitle\relax\thetitle\else\theshorttitle\fi\hfill{\pnum\folio}
\else\def\\{ and }{\pnum\folio}\hfill\ifx\theshortauthors\relax\theauthors
\else\theshortauthors\fi\fi\fi}\vss}}
\footline{\vbox to 0pt{\vglue 0mm\line{\small\pfoot\ifnum\count0=\startpage
\copyright\ \gtp\hfill\else
\agt, Volume \thevolumenumber\ (\thevolumeyear)\hfill\fi}\vss}}
\else
\font=cmsl9 scaled 1050
\font\lnum=cmbx10 
\font=cmsl9 scaled 1050
\makeatletter
\def\@oddhead{{\small\lhead\ifnum\count0=\startpage ISSN 1472-2739 
(on-line) 1472-2747 (printed)\hfill {\lnum\number\count0}\else\ifodd\count0
\def\\{ }\ifx\theshorttitle\relax \thetitle \else\theshorttitle\fi\hfill
{\lnum\number\count0}\else\def\\{ and }{\lnum\number\count0}
\hfill\ifx\theshortauthors\relax 
\theauthors\else\theshortauthors\fi\fi\fi}}\def\@evenhead{\@oddhead}
\def\@oddfoot{\small\lfoot\ifnum\count0=\startpage\copyright\ \gtp\hfill\else
\agt, Volume \thevolumenumber\ (\thevolumeyear)\hfill\fi}
\def\@evenfoot{\@oddfoot}
\makeatother
\fi
\let\maketitlepage\makeagttitle
\let\makeshorttitle\maketitlepage
\let\maketitle\maketitlepage

\newwrite\gtoutfile
\long\gdef\makeheadfile{  
{\def\\{, }\def\s{ }
\immediate\openout\gtoutfile head.xxx
\immediate\write\gtoutfile{To: math@arxiv.org}
\immediate\write\gtoutfile{Subject: put OR rep NNNNN:ppppp}
\immediate\write\gtoutfile{--text follows this line--}
\immediate\write\gtoutfile{Proxy-for: \ifx\theasciiauthors\relax
\theauthors\else\theasciiauthors\fi\s<\ifx\theasciiemail\relax\theemail\else\theasciiemail\fi>}
\immediate\write\gtoutfile{\noexpand\\}
\immediate\write\gtoutfile{Authors: \ifx\theasciiauthors\relax
\theauthors\else\theasciiauthors\fi}
{\def\\{ }\immediate\write\gtoutfile{Title: \ifx\theasciititle\relax
\thetitle\else\theasciititle\fi}}
\immediate\write\gtoutfile{Subj-class: GT or SG, GR etc}
\immediate\write\gtoutfile{MSC-class: \theprimaryclass\ifx\thesecondaryclass\relax\else, \thesecondaryclass\fi}
\immediate\write\gtoutfile{Journal-ref: Algebr. Geom. Topol. \thevolumenumber\s
(\thevolumeyear) \startpage-\finishpage}
\immediate\write\gtoutfile{Comments: Published by Algebraic and
Geometric Topology at}
\immediate\write\gtoutfile{\s\s\s  http://www.maths.warwick.ac.uk/agt/AGTVol\thevolumenumber/agt-\thevolumenumber-\thepapernumber.abs.html}
\immediate\write\gtoutfile{\noexpand\\}
\immediate\write\gtoutfile{}
\ifx\theasciiabstract\relax
\immediate\write\gtoutfile{\theabstract}\else
\immediate\write\gtoutfile{\theasciiabstract}\fi
\immediate\write\gtoutfile{}
\immediate\write\gtoutfile{\noexpand\\}
\immediate\write\gtoutfile{}
\immediate\closeout\gtoutfile}}  

\def\maketitlepage{\makeagttitle\makeheadfile}
\let\makeshorttitle\maketitlepage
\let\maketitle\maketitlepage